
\magnification=1200
\hsize=14.1cm
\vsize=19.5cm
\parindent=0cm   \parskip=0pt
\pageno=1

\def\ind{\hskip 1cm\relax}


\hoffset=15mm    
\voffset=1cm    


\ifnum\mag=\magstep1
\hoffset=-0.2cm   
\voffset=-.5cm   
\fi


\pretolerance=500 \tolerance=1000  \brokenpenalty=5000

\catcode`\@=11

\font\eightrm=cmr8         \font\eighti=cmmi8
\font\eightsy=cmsy8        \font\eightbf=cmbx8
\font\eighttt=cmtt8        \font\eightit=cmti8
\font\eightsl=cmsl8        \font\sixrm=cmr6
\font\sixi=cmmi6           \font\sixsy=cmsy6
\font\sixbf=cmbx6


\font\tengoth=eufm10       \font\tenbboard=msbm10
\font\eightgoth=eufm8      \font\eightbboard=msbm8
\font\sevengoth=eufm7      \font\sevenbboard=msbm7
\font\sixgoth=eufm6        \font\fivegoth=eufm5

\font\tencyr=wncyr10       
\font\eightcyr=wncyr8      
\font\sevencyr=wncyr7      
\font\sixcyr=wncyr6


\skewchar\eighti='177 \skewchar\sixi='177
\skewchar\eightsy='60 \skewchar\sixsy='60


\newfam\gothfam           \newfam\bboardfam
\newfam\cyrfam

\def\tenpoint{%
  \textfont0=\tenrm \scriptfont0=\sevenrm \scriptscriptfont0=\fiverm
  \def\rm{\fam\z@\tenrm}%
  \textfont1=\teni  \scriptfont1=\seveni  \scriptscriptfont1=\fivei
  \def\oldstyle{\fam\@ne\teni}\let\old=\oldstyle
  \textfont2=\tensy \scriptfont2=\sevensy \scriptscriptfont2=\fivesy
  \textfont\gothfam=\tengoth \scriptfont\gothfam=\sevengoth
  \scriptscriptfont\gothfam=\fivegoth
  \def\goth{\fam\gothfam\tengoth}%
  \textfont\bboardfam=\tenbboard \scriptfont\bboardfam=\sevenbboard
  \scriptscriptfont\bboardfam=\sevenbboard
  \def\bb{\fam\bboardfam\tenbboard}%
 \textfont\cyrfam=\tencyr \scriptfont\cyrfam=\sevencyr
  \scriptscriptfont\cyrfam=\sixcyr
  \def\cyr{\fam\cyrfam\tencyr}%
  \textfont\itfam=\tenit
  \def\it{\fam\itfam\tenit}%
  \textfont\slfam=\tensl
  \def\sl{\fam\slfam\tensl}%
  \textfont\bffam=\tenbf \scriptfont\bffam=\sevenbf
  \scriptscriptfont\bffam=\fivebf
  \def\bf{\fam\bffam\tenbf}%
  \textfont\ttfam=\tentt
  \def\tt{\fam\ttfam\tentt}%
  \abovedisplayskip=12pt plus 3pt minus 9pt
  \belowdisplayskip=\abovedisplayskip
  \abovedisplayshortskip=0pt plus 3pt
  \belowdisplayshortskip=4pt plus 3pt
  \smallskipamount=3pt plus 1pt minus 1pt
  \medskipamount=6pt plus 2pt minus 2pt
  \bigskipamount=12pt plus 4pt minus 4pt
  \normalbaselineskip=12pt
  \setbox\strutbox=\hbox{\vrule height8.5pt depth3.5pt width0pt}%
  \let\bigf@nt=\tenrm       \let\smallf@nt=\sevenrm
  \normalbaselines\rm}

\def\eightpoint{%
  \textfont0=\eightrm \scriptfont0=\sixrm \scriptscriptfont0=\fiverm
  \def\rm{\fam\z@\eightrm}%
  \textfont1=\eighti  \scriptfont1=\sixi  \scriptscriptfont1=\fivei
  \def\oldstyle{\fam\@ne\eighti}\let\old=\oldstyle
  \textfont2=\eightsy \scriptfont2=\sixsy \scriptscriptfont2=\fivesy
  \textfont\gothfam=\eightgoth \scriptfont\gothfam=\sixgoth
  \scriptscriptfont\gothfam=\fivegoth
  \def\goth{\fam\gothfam\eightgoth}%
  \textfont\cyrfam=\eightcyr \scriptfont\cyrfam=\sixcyr
  \scriptscriptfont\cyrfam=\sixcyr
  \def\cyr{\fam\cyrfam\eightcyr}%
  \textfont\bboardfam=\eightbboard \scriptfont\bboardfam=\sevenbboard
  \scriptscriptfont\bboardfam=\sevenbboard
  \def\bb{\fam\bboardfam}%
  \textfont\itfam=\eightit
  \def\it{\fam\itfam\eightit}%
  \textfont\slfam=\eightsl
  \def\sl{\fam\slfam\eightsl}%
  \textfont\bffam=\eightbf \scriptfont\bffam=\sixbf
  \scriptscriptfont\bffam=\fivebf
  \def\bf{\fam\bffam\eightbf}%
  \textfont\ttfam=\eighttt
  \def\tt{\fam\ttfam\eighttt}%
  \abovedisplayskip=9pt plus 3pt minus 9pt
  \belowdisplayskip=\abovedisplayskip
  \abovedisplayshortskip=0pt plus 3pt
  \belowdisplayshortskip=3pt plus 3pt
  \smallskipamount=2pt plus 1pt minus 1pt
  \medskipamount=4pt plus 2pt minus 1pt
  \bigskipamount=9pt plus 3pt minus 3pt
  \normalbaselineskip=9pt
  \setbox\strutbox=\hbox{\vrule height7pt depth2pt width0pt}%
  \let\bigf@nt=\eightrm     \let\smallf@nt=\sixrm
  \normalbaselines\rm}

\tenpoint


\def\pc#1{\bigf@nt#1\smallf@nt}         \def\pd#1 {{\pc#1} }


\catcode`\;=\active
\def;{\relax\ifhmode\ifdim\lastskip>\z@\unskip\fi
\kern\fontdimen2  -1.2 \fontdimen3 \string;}

\catcode`\:=\active
\def:{\relax\ifhmode\ifdim\lastskip>\z@\unskip\fi\penalty\@M\ \fi\string:}

\catcode`\!=\active
\def!{\relax\ifhmode\ifdim\lastskip>\z@
\unskip\fi\kern\fontdimen2  -1.1 \fontdimen3 \string!}

\catcode`\?=\active
\def?{\relax\ifhmode\ifdim\lastskip>\z@
\unskip\fi\kern\fontdimen2  -1.1 \fontdimen3 \string?}

\def\^#1{\if#1i{\accent"5E\i}\else{\accent"5E #1}\fi}
\def\"#1{\if#1i{\accent"7F\i}\else{\accent"7F #1}\fi}

\frenchspacing


\newtoks\auteurcourant      \auteurcourant={\hfil}
\newtoks\titrecourant       \titrecourant={\hfil}

\newtoks\hautpagetitre      \hautpagetitre={\hfil}
\newtoks\baspagetitre       \baspagetitre={\hfil}

\newtoks\hautpagegauche
\hautpagegauche={\eightpoint\rlap{\folio}\hfil\the\auteurcourant\hfil}
\newtoks\hautpagedroite
\hautpagedroite={\eightpoint\hfil\the\titrecourant\hfil\llap{\folio}}

\newtoks\baspagegauche      \baspagegauche={\hfil}
\newtoks\baspagedroite      \baspagedroite={\hfil}

\newif\ifpagetitre          \pagetitretrue

a marche. Alors ...


\footline={\ifpagetitre\the\baspagetitre\else
\ifodd\pageno\the\baspagedroite\else\the\baspagegauche\fi\fi
\global\pagetitrefalse}


\def\raggedbottom{\topskip 10pt plus 36pt\r@ggedbottomtrue}



\def\pointir{\unskip . --- \ignorespaces}


\def\Bigbreak{\vskip-\lastskip\bigbreak}
\def\Medbreak{\vskip-\lastskip\medbreak}


\def\ctexte#1\endctexte{%
  \hbox{$\vcenter{\halign{\hfill##\hfill\crcr#1\crcr}}$}}


\long\def\ctitre#1\endctitre{%
    \ifdim\lastskip<24pt\vskip-\lastskip\bigbreak\bigbreak\fi
  		\vbox{\parindent=0pt\leftskip=0pt plus 1fill
          \rightskip=\leftskip
          \parfillskip=0pt\bf#1\par}
    \bigskip\nobreak}

\long\def\section#1\endsection{%
\vskip 0pt plus 3\normalbaselineskip
\penalty-250
\vskip 0pt plus -3\normalbaselineskip
\Bigbreak
\message{[section \string: #1]}{\bf#1\unskip}\pointir}

\long\def\sectiona#1\endsection{%
\vskip 0pt plus 3\normalbaselineskip
\penalty-250
\vskip 0pt plus -3\normalbaselineskip
\Bigbreak
\message{[sectiona \string: #1]}%
{\bf#1}\medskip\nobreak}

\long\def\subsection#1\endsubsection{%
\Medbreak
{\it#1\unskip}\pointir}

\long\def\subsectiona#1\endsubsection{%
\Medbreak
{\it#1}\par\nobreak}

\def\rem#1\endrem{%
\Medbreak
{\it#1\unskip} : }

\def\remp#1\endrem{%
\Medbreak
{\pc #1\unskip}\pointir}

\def\rema#1\endrem{%
\Medbreak
{\it #1}\par\nobreak}

\def\newparwithcolon#1\endnewparwithcolon{
\Medbreak
{#1\unskip} : }

\def\newparwithpointir#1\endnewparwithpointir{
\Medbreak
{#1\unskip}\pointir}

\def\newpara#1\endnewpar{
\Medbreak
{#1\unskip}\smallskip\nobreak}

\let\+=\tabalign

\def\signature#1\endsignature{\vskip 15mm minus 5mm\rightline{\vtop{#1}}}

\mathcode`A="7041 \mathcode`B="7042 \mathcode`C="7043 \mathcode`D="7044
\mathcode`E="7045 \mathcode`F="7046 \mathcode`G="7047 \mathcode`H="7048
\mathcode`I="7049 \mathcode`J="704A \mathcode`K="704B \mathcode`L="704C
\mathcode`M="704D \mathcode`N="704E \mathcode`O="704F \mathcode`P="7050
\mathcode`Q="7051 \mathcode`R="7052 \mathcode`S="7053 \mathcode`T="7054
\mathcode`U="7055 \mathcode`V="7056 \mathcode`W="7057 \mathcode`X="7058
\mathcode`Y="7059 \mathcode`Z="705A

\def\spacedmath#1{\def\packedmath##1${\bgroup\mathsurround=0pt ##1\egroup$}%
\mathsurround#1 \everymath={\packedmath}\everydisplay={\mathsurround=0pt }}

\def\nospacedmath{\mathsurround=0pt \everymath={}\everydisplay={} }


\long\def\th#1 #2\enonce#3\endth{%
   \Medbreak
   {\pc#1} {#2\unskip}\pointir{\it #3}\medskip}

\long\def\tha#1 #2\enonce#3\endth{%
   \Medbreak
   {\pc#1} {#2\unskip}\par\nobreak{\it #3}\medskip}


\long\def\Th#1 #2 #3\enonce#4\endth{%
   \Medbreak
   #1 {\pc#2} {#3\unskip}\pointir{\it #4}\medskip}

\long\def\Tha#1 #2 #3\enonce#4\endth{%
   \Medbreak
   #1 {\pc#2} #3\par\nobreak{\it #4}\medskip}


\def\decale#1{\smallbreak\hskip 28pt\llap{#1}\kern 5pt}
\def\decaledecale#1{\smallbreak\hskip 34pt\llap{#1}\kern 5pt}
\def\puce{\smallbreak\hskip 6pt{$\scriptstyle\bullet$}\kern 5pt}



\def\displaylinesno#1{\displ@y\halign{
\hbox to\displaywidth{$\@lign\hfil\displaystyle##\hfil$}&
\llap{$##$}\crcr#1\crcr}}


\def\ldisplaylinesno#1{\displ@y\halign{
\hbox to\displaywidth{$\@lign\hfil\displaystyle##\hfil$}&
\kern-\displaywidth\rlap{$##$}\tabskip\displaywidth\crcr#1\crcr}}


\def\eqalign#1{\null\,\vcenter{\openup\jot\m@th\ialign{
\strut\hfil$\displaystyle{##}$&$\displaystyle{{}##}$\hfil
&&\quad\strut\hfil$\displaystyle{##}$&$\displaystyle{{}##}$\hfil
\crcr#1\crcr}}\,}

\def\system#1{\left\{\null\,\vcenter{\openup1\jot\m@th
\ialign{\strut$##$&\hfil$##$&$##$\hfil&&
        \enskip$##$\enskip&\hfil$##$&$##$\hfil\crcr#1\crcr}}\right.}


\let\@ldmessage=\message

\def\message#1{{\def\pc{\string\pc\space}%
                \def\'{\string'}\def\`{\string`}%
                \def\^{\string^}\def\"{\string"}%
                \@ldmessage{#1}}}

\def\diagram#1{\def\normalbaselines{\baselineskip=0pt
\lineskip=10pt\lineskiplimit=1pt}   \matrix{#1}}



\def\up#1{\raise 1ex\hbox{\smallf@nt#1}}


\def\cf{{\it cf}} 

\def\qed{\raise -2pt\hbox{\vrule\vbox to 10pt{\hrule width 4pt
                 \vfill\hrule}\vrule}}

\def\cqfd{\unskip\penalty 500\quad\vrule height 4pt depth 0pt width
4pt\medbreak}

\def\virg{\raise .4ex\hbox{,}}   


\def\build#1_#2^#3{\mathrel{
\mathop{\kern 0pt#1}\limits_{#2}^{#3}}}


\def\boxit#1#2{%
\setbox1=\hbox{\kern#1{#2}\kern#1}%
\dimen1=\ht1 \advance\dimen1 by #1 \dimen2=\dp1 \advance\dimen2 by #1
\setbox1=\hbox{\vrule height\dimen1 depth\dimen2\box1\vrule}%
\setbox1=\vbox{\hrule\box1\hrule}%
\advance\dimen1 by .4pt \ht1=\dimen1
\advance\dimen2 by .4pt \dp1=\dimen2  \box1\relax}

\def\date{\the\day\ \ifcase\month\or janvier\or f\'evrier\or mars\or
avril\or mai\or juin\or juillet\or ao\^ut\or septembre\or octobre\or
novembre\or d\'ecembre\fi \ {\old \the\year}}

\def\dateam{\ifcase\month\or January\or February\or March\or
April\or May\or June\or July\or August\or September\or October\or
November\or December\fi \ \the\day ,\ \the\year}

\def\crog{{\vrule height 2.57mm depth 0.85mm width 0.3mm}\kern -0.36mm
[}

\def\crod{]\kern -0.4mm{\vrule height 2.57mm depth 0.85mm
width 0.3 mm}}

\def\rond{\kern 1pt{\scriptstyle\circ}\kern 1pt}

\def\diagram#1{\def\normalbaselines{\baselineskip=0pt
\lineskip=10pt\lineskiplimit=1pt}   \matrix{#1}}

\def\hfl#1#2{\nospacedmath\smash{\mathop{\hbox to
12mm{\rightarrowfill}}\limits^{\scriptstyle#1}_{\scriptstyle#2}}}

\def\ghfl#1#2{\nospacedmath\smash{\mathop{\hbox to
25mm{\rightarrowfill}}\limits^{\scriptstyle#1}_{\scriptstyle#2}}}

\def\phfl#1#2{\nospacedmath\smash{\mathop{\hbox to
8mm{\rightarrowfill}}\limits^{\scriptstyle#1}_{\scriptstyle#2}}}

\def\vfl#1#2{\llap{$\scriptstyle#1$}\left\downarrow\vbox to
6mm{}\right.\rlap{$\scriptstyle#2$}}

\def\va{vari\'et\'e ab\'elienne}
\def\vas{vari\'et\'es  ab\'eliennes}

\def\pa{\S\kern.15em}
\def\ra{\rightarrow}
\def\lra{\longrightarrow}
\def\llra{\nospacedmath\hbox to 10mm{\rightarrowfill}}
\def\lllra{\nospacedmath\hbox to 15mm{\rightarrowfill}}
\def\saut{\vskip 5mm plus 1mm minus 2mm}
\def\Z{\hbox{\bf Z}}
\def\P{\hbox{\bf P}}

\def\C{\hbox{\bf C}}

\def\isom{\simeq}

\def\Tr{\mathop{\rm Tr}\nolimits}

\def\Hom{\mathop{\rm Hom}\nolimits}

\def\Sing{\mathop{\rm Sing}\nolimits}
\def\Ker{\mathop{\rm Ker}\nolimits}
\def\Card{\mathop{\rm Card}\nolimits}

\def\Pic{\mathop{\rm Pic}\nolimits}

\def\codim{\mathop{\rm codim}\nolimits}

\def\loc{{\it loc.cit.\/}}
\def\cf{{\it cf.\/}}
\def\cad{c'est-\`a-dire}
\def\ssi{si et seulement si}

\def\cc#1{\hfill\kern .7em#1\kern .7em\hfill}

\def\theo{th\'eor\`eme}
\def\limproj{\mathop{\oalign{lim\cr\hidewidth$\longleftarrow$\hidewidth\cr}}}
\def\dra{\ra\kern -3mm\ra}
\def\ldra{\lra\kern -3mm\ra}


\catcode`\@=12

\showboxbreadth=-1  \showboxdepth=-1



\baselineskip=14pt
\spacedmath{1.7pt}
\baspagegauche={\centerline{\tenbf\folio}}
\baspagedroite={\centerline{\tenbf\folio}}
\hautpagegauche={\hfil}
\hautpagedroite={\hfil}
\font\eightrm=cmr10 at 8pt
\parskip=1.7mm

\def\cc#1{\hfill\kern .7em#1\kern .7em\hfill}
\def\theo{th\'eor\`eme}
\def\pp{\pi_1^{\rm mod}}

\saut
\ctitre
{\bf THEOREMES DE CONNEXITE ET VARIETES ABELIENNES}
\endctitre
\medskip
\centerline{{Olivier {\pc DEBARRE}}
\footnote{(*)}{\rm Financ\'e en partie par
N.S.F. Grant DMS 92-03919 et le Projet Europ\'een Science
 ``Geometry of Algebraic Varieties", Contract no. SCI-0398-C (A).}}

\vskip 1cm
\ind Cet article est consacr\'e \`a la d\'emonstration d'un \theo\
de connexit\'e pour les \vas , ainsi qu'\`a l'exposition de quelques
corollaires. Les r\'esultats analogues pour les espaces projectifs et leurs
cons\'equences sont connus depuis un moment d\'ej\`a ([FH], [FL], [F1]) et
fournissent un guide qu'il n'y a qu'\`a suivre.

\ind Le \theo\ principal \'enonce que si on se donne une \va\ $X$ et
deux morphismes $V\ra X$ et $W\ra X$, alors, sous certaines hypoth\`eses, assez
contraignantes mais indispensables, $V\times _XW$ est {\it connexe}. Les
cons\'equences de ce
r\'esultat sont nombreuses.

\ind Par exemple, on montre que le groupe fondamental alg\'ebrique d'une
sous-vari\'et\'e {\it normale}
d'une \va\ {\it simple} $X$, de dimension $>{1\over 2}\dim (X)$, est isomorphe
\`a celui de $X$
(cela d\'ecoule de r\'esultats ant\'erieurs ([S]) lorsque $V$ est lisse).

\ind On \'etend aussi certains r\'esultats de Nori ([N]): si $D$ est un
diviseur dans une
\va\ $X$, dont aucune composante n'est une sous-\va , et qui est \`a
croisements normaux en
dehors d'un ferm\'e de codimension $2$ dans $D$, alors le noyau du morphisme
surjectif $\pi_1(X-D)\ra
\pi_1(X)$ est ab\'elien libre de type fini, et l'on d\'etermine son rang.

\ind Des r\'esultats analogues \`a ceux de [L1] sont aussi obtenus, comme par
exemple que le groupe
fondamental alg\'ebrique d'une vari\'et\'e {\it normale} $V$ qui est
rev\^etement de degr\'e $\le \dim
(V)$ d'une \va\ simple
$X$, est isomorphe \`a celui de $X$. Lorsque $\dim (V)\ge 2$, la m\^eme
conclusion subsiste pour un
rev\^etement qui induit une bijection ensembliste au-dessus d'une courbe, comme
par exemple
un rev\^etement cyclique, sans condition sur son degr\'e.

\ind On termine avec une conjecture, analogue d'un r\'esultat d\'emontr\'e pour
les espaces
projectifs par Lazarsfeld. Nous renvoyons le lecteur au \pa 8 pour
son \'enonc\'e, un peu technique.

{\bf 1. Rappels et notations}

\ind Dans tout cet article, on travaille sur un corps alg\'ebriquement clos $k$
{\it de
caract\'eristique nulle}. Une {\it vari\'et\'e} est un sch\'ema projectif
r\'eduit de type fini sur
$k$, pas n\'ecessairement irr\'eductible, mais qu'on supposera toujours non
vide.

\ind Le groupe fondamental alg\'ebrique d'un sch\'ema est d\'efini dans [G1],
Exp. V (\cf\ aussi [Mu1],
p. 169). Comme nous ne consid\'ererons que des espaces connexes, nous
\'ecrirons simplement $\pi_1^{\rm
alg}(X)$ pour le groupe fondamental alg\'ebrique d'une vari\'et\'e connexe $X$.

\ind Lorsque $k=\C $, le groupe $\pi_1^{\rm alg}(X)$ est le compl\'et\'e du
groupe fondamental
topologique $\pi_1(X)$ pour la topologie des sous-groupes d'indices finis
([G1], Exp. XII, cor.
5.2).

(1.1) Soit $\Pic (X)$ le sch\'ema de Picard de $X$. Pour tout
entier $n>0$, on note ${}_n\Pic (X)$ le noyau de la multiplication par $n$ dans
$\Pic (X)$; il existe
un isomorphisme ([G1], Exp. XI, p. 305): $$\Hom ( \pi_1^{\rm alg}(X),\Z /n\Z
)\isom {}_n\Pic (X)\ .$$

\ind On notera $\Pic ^0(X)$ la composante neutre de $\Pic (X)$. Lorsque $X$ est
normal, $\Pic
^0(X)_{\rm red}$ est une vari\'et\'e ab\'elienne ([G2], Exp. VI, cor. 3.2).

\ind Soient $V$ une vari\'et\'e connexe et $f:V\ra X$ un morphisme. Il induit
alors un morphisme
$f_*:\pi_1^{\rm alg}(V)\ra\pi_1^{\rm alg}(X)$ ([G1], Exp. V, p. 142) et:

(1.2) $f_*$ est surjectif \ssi\ pour tout rev\^etement \'etale connexe $X'\ra
X$, le rev\^etement
\'etale $V\times _XX'\ra V$ est connexe ([G1], Exp. V, prop. 6.9);

(1.3) $f_*$ est injectif \ssi\ pour tout rev\^etement \'etale $V'\ra V$, il
existe un rev\^etement
\'etale $X'\ra X$ et un $V$--morphisme d'une composante connexe de $V\times
_XX'$ dans $V'$. C'est
le cas en particulier si tout rev\^etement \'etale connexe de $V$ est isomorphe
\`a un rev\^etement
du type $V\times _XX'\ra V$ ([G1], Exp. V, cor. 6.8 et remarque 6.12);

(1.4) si $f$ est \'etale, alors $f_*$ est injectif ([G1], Exp. V, p. 142).

(1.5) On suppose de plus $X$ normal. Lorsque $f_*$ est {\it surjectif}, le
morphisme induit $f^*_0:\Pic ^0(X)_{\rm red}\ra\Pic ^0(V)_{\rm red}$ est un
{\it plongement}. En
effet, pour tout entier $n>0$, il ressort de (1.1) que les morphismes induits
$f^*_n:{}_n\Pic
(X)\ra{}_n\Pic (V)$ sont surjectifs, de sorte que $\Ker (f^*_0)$ est sans
torsion, donc nul. Si de
plus $V$ est normal et $f_*$ bijectif, alors $f^*_0$ est un isomorphisme
(m\^eme d\'emonstration).

(1.6) Supposons $X$ {\it (g\'eom\'etriquement) unibranche}, \cad\ que pour tout
point $x$ de $X$, le
normalis\'e de l'anneau local ${\cal O}_{X,x}$ a un seul id\'eal maximal.
Alors, si $X'\ra X$ est un
rev\^etement \'etale connexe, $X'$ est irr\'eductible et unibranche ([GD1],
17.5.7). D'autre part,
si $\tilde X\ra X$ est la normalisation; le morphisme induit $\pi_1^{\rm
alg}(\tilde X)\ra\pi_1^{\rm
alg}(X)$ est un isomorphisme ([G1], Exp. IX, th. 4.10).

(1.7) Soit maintenant $X$ une \va . Alors $\pi_1^{\rm alg}(X)$ est isomorphe
\`a $\hat{\Z }^{2\dim
(X)}$ ([G1], Exp. XI, th. 2.1). Il r\'esulte de [Mu1], p. 147 et 171, que deux
\vas\ isog\`enes
ont des groupes fondamentaux alg\'ebriques isomorphes. Lorsque $k=\C $, le
groupe $\pi_1(X)$ est
isomorphe \`a $\Z ^{2\dim (X)}$.

(1.8) Soient $V$ une vari\'et\'e et $f:V\ra X$ un morphisme. On appelera
factorisation de $f$ la
donn\'ee d'une \va\ $X'$, d'une isog\'enie $q:X'\ra X$ et d'un morphisme
$f':V\ra X'$
tels que $f=qf'$. On dira que $f$ est {\it minimal} si, pour toute
factorisation, $q$ est un
isomorphisme. Supposons que $f(V)$
engendre $X$. Le noyau du morphisme induit $f^*:\Pic ^0(X)\ra \Pic (V)$ est
alors fini et correspond
\`a une isog\'enie $q:X'\ra X$ \`a travers laquelle $f$ se factorise. La
factorisation correspondante
de $f$ est alors minimale.

\ind Enfin, si $V$ est unibranche et $f:V\ra X$ minimal, et si $\eta:\tilde
V\ra V$ est
la normalisation de $V$, le morphisme compos\'e $\tilde f=f\eta$ est aussi
minimal.
En effet, soit $\tilde V\ra X'\ra X$ une factorisation de $\tilde f$. Le
morphisme $\eta$ se factorise
alors en $\tilde V\buildrel{g}\over{\ra}V\times _XX'\ra V$.  Soit $V'$ la
composante connexe de
$V\times _XX'$ qui contient $g(\tilde V)$. Par (1.6), $V'$ est irr\'eductible,
donc \'egal \`a
$g(\tilde V)$. Le morphisme \'etale $V'\ra V$ est donc un isomorphisme; $f$ se
factorise alors \`a
travers $q$, qui est donc un isomorphisme.

(1.9) Pour \'enoncer notre r\'esultat principal (\theo\ 4.5) avec des
hypoth\`eses
minimales, nous aurons besoin de la d\'efinition suivante. Soit $X$ une \va\ et
soient $V$ et $W$
deux sous-vari\'et\'es de $X$. On dit que $(V,W)$ {\it remplit} (resp. {\it
remplit strictement}) $X$
si, pour toute sous-\va\ propre $K$ de $X$, on a: $$\dim \bigl( \pi (V)\bigr)
+\dim \bigl( \pi
(W)\bigr)\ge \dim (X/K)$$
(resp. $>$), o\`u $\pi :X\ra X/K$ est la surjection canonique. Par exemple,
$(V,X)$  remplit
strictement $X$ \ssi\ $V$ engendre $X$. Lorsque $X$ est simple, ou encore
lorsque toute courbe
contenue dans $V$ engendre $X$, le couple $(V,W)$ remplit (resp. remplit
strictement) $X$ \ssi\ $\dim
(V) +\dim (W)\ge \dim (X)$ (resp. $>$).

\ind Afin de ne pas trop alourdir le texte, nous avons pr\'ef\'er\'e \'enoncer
les diff\'erents corollaires avec des hypoth\`eses plus faibles (mais plus
parlantes). Le lecteur
pourra r\'etablir facilement les \'enonc\'es optimaux.

{\bf 2. Sous-vari\'et\'es lisses}

\ind En ce qui concerne les sous-vari\'et\'es {\it lisses} d'une \va ,
Sommese a obtenu dans une s\'erie d'articles, qui culmine avec [So], des
r\'esultats tr\`es complets
dont voici un \'echantillon:

{\bf Th\'eor\`eme 2.1.} (Sommese) -- {\it Soit $X$ une \va\ complexe simple et
soient $V$ et $W$ deux sous-vari\'et\'es lisses de $X$. On suppose $W$ connexe
et $V\cap W$ non vide.
Alors  $\pi_q(W,V\cap W)=0$ pour $q\le \min (\dim (W),\dim (V)+1)
-\codim (V)$.}

\ind On verra plus loin (\theo\ 3.1) que $V\cap W$ est non vide d\`es que $\dim
(W)\ge \codim
(V)$, \cad\ d\`es que la borne sup\'erieure sur $q$ est $\ge 0$.

{\bf Corollaire 2.2.} -- {\it Soient $X$ une \va\ complexe simple et $V$ une
sous-vari\'et\'e lisse connexe de $X$. Alors:
{\parindent=1cm
\item{\rm (i)} le morphisme
induit $\pi_1(V)\ra\pi_1(X)$ est surjectif pour
 $\dim (V)\ge{1\over 2}\dim (X)$ et est un isomorphisme pour  $\dim (V)>{1\over
2}\dim (X)$,
\item{\rm (ii)} le morphisme de restriction $\Pic (X)\ra\Pic (V)$ est injectif
pour
 $\dim (V)>{1\over 2}\dim (X)$ et est un isomorphisme pour  $\dim (V)\ge
1+{1\over 2}\dim (X)$.\par}}

{\bf Remarque 2.3.} Soit $C$ une courbe projective
lisse connexe de genre $g$, soit $JC$ sa jacobienne et soit $W_d(C)$ la
sous-vari\'et\'e de $JC$
qui param\`etre les classes d'\'equivalence de diviseurs effectifs de degr\'e
$d$ sur $C$.
Lorsque $d<g$, pour  tout point $x$ de $C$, le diviseur de Weil $x+W_{d-1}(C)$
ne peut \^etre la
restriction \`a $W_d(C)$ d'un diviseur de $JC$, pour des raisons num\'eriques.
Pour $C$
g\'en\'erique et $d<1+{1\over 2}g$, $JC$ est simple, $W_d(C)$ est lisse, mais
la restriction $\Pic
(JC)\ra\Pic\bigl( W_d(C)\bigr)$ n'est donc pas surjective. La deuxi\`eme borne
du corollaire 2.2.(ii)
est donc la meilleure possible. On remarquera qu'on retrouve aussi ainsi le
fait bien connu que
$W_d(C)$ n'est pas lisse pour $g>d\ge 1+{1\over 2}g$. En g\'en\'eral, on peut
montrer que la
restriction $\Pic (JC)\ra\Pic\bigl( W_d(C)\bigr)$ est un isomorphisme \ssi\
$W_d(C)$ est singulier,
donc en particulier lorsque $d\ge 1+{1\over 2}g$. Je pense d'ailleurs que le
(ii) du
corollaire devrait rester vrai pour $V$ normal.

{\bf 3. Sous-vari\'et\'es quelconques}

\ind Les premiers r\'esultats sur l'intersection de deux sous-vari\'et\'es
quelconques d'une \va\
pr\'ec\`edent en fait ceux de Sommese. Ils sont dus \`a Barth, qui a
d\'emontr\'e dans [B1] le
\theo\ suivant lorsque $X$ est simple.

{\bf Th\'eor\`eme 3.1.} -- {\it Soit $X$ une \va\ et soient $V$ et $W$ deux
sous-vari\'et\'es irr\'eductibles de $X$ telles que
$(V,W)$ remplisse $X$ (\cf\ (1.9)). Alors $V\cap W\ne \emptyset$.}

{\bf Remarques 3.2.} 1) La conclusion ne subsiste pas en g\'en\'eral si $V$ ou
$W$ est
r\'eductible: si $Y$ est une \va\  simple de dimension $\ge 2$, $E$ une courbe
elliptique, $C$ une
courbe dans $Y$ ne contenant pas l'origine, et $e$ un point non nul de
$E$, alors $V=C\times\{ 0\}$ et $W=\{ 0\}\times E\cup Y\times \{ e\}$
remplissent $X=Y\times E$, mais
ne se rencontrent pas.

\ind 2) Le d\'ebut de la d\'emonstration montre que lorsque toute courbe
contenue dans $V$ engendre $X$, la conclusion subsiste
quelle que soit la caract\'eristique de $k$. De plus, si $V_1,\ldots ,V_r$ sont
des
sous-vari\'et\'es de $X$ v\'erifiant $\sum_{i=1}^r\codim _X(V_i)\le \dim (X)$
et si toute courbe
contenue dans $V_1$ engendre $X$, alors
$\bigcap_{i=1}^rV_i\ne \emptyset$.

{\bf D\'emonstration du \theo .}  On raisonne par r\'ecurrence sur la dimension
de $X$. Si
$V$ ne rencontre pas $W$, le point $0$ n'est pas dans l'image $(V-W)$ du
morphisme $m:V\times W\ra X$
d\'efini par $m(v,w)=v-w$. Ce dernier n'est donc pas surjectif. Soit $a=v-w$ un
point
g\'en\'erique, donc lisse, de  $(V-W)$. Il existe une sous-vari\'et\'e $F$ de
$X$ de dimension
$>0$, contenant $0$, telle que:
$$m^{-1}(a)=\{\ (v+e,w+e)\bigm| e\in F\ \}\ .$$
\ind On a alors $(a+F-F)\subset (V-W)$, de sorte que $T_a(a+F-F)\subset
T_{a}(V-W)$.
Soit $\! <\! F\! >\! $ la sous-\va\ de $X$ engendr\'ee par $F$.

{\bf Lemme 3.3.} -- {\it On a $T_0(F-F)=T_0\! <\! F\! >\! $.}

{\bf D\'emonstration.} En utilisant le \theo\ de Lefschetz, on se ram\`ene \`a
montrer que si une
courbe $C$ (peut-\^etre r\'eductible) engendre $X$, alors $T_0(C-C)=T_0X$. On
note\break ${\cal
G}:C_{\rm reg}\ra\P T_0X$ l'application de Gauss. Pour tout $x\in C_{\rm reg}$,
on a:
$${\cal G}(x)=\P T_0(C-x)\subset\P T_0(C-C)\ ,$$de sorte qu'il suffit de
montrer que l'image de
${\cal G}$ n'est pas contenue dans un hyperplan. Cela r\'esulte du fait que la
restriction
$H^0(X,\Omega^1_X)\ra H^0(C_{\rm reg},\Omega^1_{C_{\rm reg}})$ est
injective.\cqfd

\ind On a donc $T_a(a+\! <\! F\! >\! )\subset T_a(V-W)$. Le point $a$
\'etant g\'en\'erique, on se convainc rapidement que pour $x$ dans un voisinage
de $a$ dans $(V-W)$,
la sous-\va\ correspondante (\cad\ $\! <\! pr_1m^{-1}(x)\! >\! $) reste \'egale
\`a $\! <\! F\! >\!
$. Pour $x$ g\'en\'erique dans $(V-W)$, on a donc $T_x(x+\! <\! F\! >\!
)\subset T_x(V-W)$.

{\bf Lemme 3.4.} -- {\it Soient $X$ une \va , $Z$ une sous-vari\'et\'e
irr\'eductible de $X$ et $K$
une sous-\va\ de $X$. On suppose que, pour $z$ g\'en\'erique dans $Z$, on a
$T_z(z+K)\subset T_zZ$.
Alors, $Z+K=Z$.}

{\bf D\'emonstration.} Comme on est en caract\'eristique $0$, il existe un
point
$(z,k)$ lisse sur $Z\times K$ en lequel la diff\'erentielle du morphisme
d'addition $Z\times K\ra (Z+K)$ est surjective. L'image de cette
diff\'erentielle est
par hypoth\`ese $T_{z+k}(Z+k)$, de sorte que $\dim (Z+K)=\dim (Z)$, ce qui
prouve le
lemme.\cqfd

\ind Il est clair que les images de $V$ et $W$ dans $X/\! <\! F\! >\! $
remplissent $X/\! <\! F\! >\! $. Le lemme entra\^ine alors $V-W+\! <\! F\! >\,
=V-W$ et l'hypoth\`ese de
r\'ecurrence $(V-W)\,\cap
\! <\! F\! >\, \neq \emptyset$. On en d\'eduit $0\in (V-W)$, ce qui termine la
d\'emonstration.\cqfd

{\bf Corollaire 3.5.} -- {\it Soient $X$ une \va\ simple et $V$ une
sous-vari\'et\'e de $X$
v\'erifiant $\dim (V)\ge {1\over 2}\dim (X)$. Alors le morphisme induit
$\pi_1^{\rm
alg}(V)\ra\pi_1^{\rm alg}(X)$ est surjectif. Lorsque $k=\C$, il en est de
m\^eme du morphisme
$\pi_1(V)\ra\pi_1(X)$.}

{\bf D\'emonstration.} Soit $X'\ra X$ une isog\'enie. Par (1.2), il s'agit de
montrer que
la sous-vari\'et\'e $V\times _XX'$ de $X'$ est connexe. Cela r\'esulte du fait
que ses composantes
irr\'eductibles sont toutes de m\^eme dimension que $V$, donc se rencontrent
deux \`a deux par le
\theo . Lorsque $k=\C$, soit $\Gamma$ le $\Z$--module libre $\pi_1(X)$ (1.7).
Pour tout entier $n>0$,
le morphisme compos\'e $\pi_1(V)\ra\Gamma\ra\Gamma/n\Gamma$ est alors
surjectif. Le conoyau $Q$ de
$\pi_1(V)\ra\Gamma$ est alors un groupe ab\'elien de type fini qui v\'erifie
$Q=nQ$ pour tout $n>0$.
Il est donc nul, ce qui termine la d\'emonstration du corollaire.\cqfd

{\bf Remarques 3.6.} 1) L'hypoth\`ese optimale est ici que $(V,V)$ remplit $X$
(\cf\ (1.9)).

\ind 2) Sous les hypoth\`eses du corollaire, on a en particulier que le
morphisme d'inclusion
$V\hookrightarrow X$ est minimal au sens de (1.8).

\ind Dans un second article [B2], semble-t-il peu connu, Barth s'est ensuite
int\'eress\'e \`a la
connexit\'e de l'intersection de deux sous-vari\'et\'es, soit dans une \va ,
soit dans l'espace
projectif, obtenant ainsi des r\'esultats dont il faudra attendre plus de dix
ans pour
qu'ils soient g\'en\'eralis\'es et appliqu\'es par Fulton et Hansen dans le cas
de l'espace
projectif (\cf\ [F1] et [FL] pour des compte-rendus). A ma connaissance, le cas
des \vas\ est
rest\'e ouvert. Il s'av\`ere que beaucoup des r\'esultats obtenus dans l'espace
projectif restent
valables avec quelques modifications.

{\bf 4. Th\'eor\`emes de connexit\'e}

\ind Commen\c cons par le r\'esultat technique suivant, dont la d\'emonstration
suit des id\'ees de
Mumford telles qu'elles sont expos\'ees dans [FL], p. 42.

{\bf Proposition 4.1.} -- {\it Soient $X$ une \va , $\tilde V$, $\tilde W$ et
$Y$ des vari\'et\'es
irr\'eductibles normales, $f:\tilde V\ra X$ et $g:\tilde W\ra X$ deux
morphismes et $\pi:Y\ra \tilde V\times\tilde W$ un morphisme \'etale. On
suppose que $\bigl( f(\tilde
V),g(\tilde W)\bigr)$ remplit strictement $X$. Soit $m:X\times X\ra X$ le
morphisme d\'efini par
$m(x,y)=x-y$. Alors le morphisme $M=m\circ (f,g)\circ\pi$ se factorise en
$M:Y\buildrel{h}\over{\ra}
X'\buildrel{q}\over{\ra} X$, o\`u $X'$ est une \va , $q$ une isog\'enie, et
o\`u les
fibres de $h$ sont connexes.}

{\bf D\'emonstration.} Elle est bas\'ee sur le lemme suivant:

{\bf Lemme 4.2} -- {\it Sous les hypoth\`eses de la proposition, si $D$ est une
sous-vari\'et\'e
irr\'eductible de $X$ et $Z$ une composante irr\'eductible de $M^{-1}(D)$ telle
que, pour $z$
g\'en\'erique dans $Z$, on ait $T_zM(T_zZ)\subset T_{M(z)}D$, alors $D$ n'est
pas un diviseur.}

{\bf D\'emonstration.} Proc\'edant par r\'ecurrence sur la dimension de
$X$, on raisonne par l'absurde en supposant que $D$ est un diviseur. Posons
$V'=f(\tilde V)$ et
$W'=g(\tilde W)$. On note encore $m$ l'application $V'\times W'\ra X$ induite
par $m$. Soit
$z$ un point g\'en\'erique de $Z$, on note $\pi (z)=(\tilde v,\tilde w)$,
$v'=f(\tilde v)$,
$w'=g(\tilde w)$ et $a=M(z)=v'-w'$. Comme $\dim (V')+\dim (W')>\dim (X)$, il
existe une
sous-vari\'et\'e $F$ de $X$ de dimension $>0$, contenant $0$, telle que:
$$m^{-1}(a)=\{\
(v'+e,w'+e)\bigm| e\in F\ \}\ .$$ \ind Le lieu o\`u l'application
$\tilde V\ra V'$ induite par $f$ n'est pas lisse, est un ferm\'e propre de
$\tilde V$. Il en est de
m\^eme pour le morphisme $\tilde W\ra W'$ induit par $g$. Comme $Z$ est un
diviseur et que $T$ est
\'etale, on peut donc supposer par exemple que l'image de la diff\'erentielle
de $(f,g)\circ\pi$ en
$z$ contient $T_{v'}(V')\times\{ 0\}$. Il en r\'esulte que l'on a, pour tout
point $e$ de $F$:
$$a+T_eF\subset -w'+T_{v'}V'\subset T_zM(T_zZ)\subset T_aD\ .$$ \ind Si $K$ est
la sous-\va\ de $X$
(de dimension $>0$) engendr\'ee par $F$, il en r\'esulte que $T_a(a+K)\subset
T_aD$. Comme dans la
d\'emonstration du \theo\ 3.1, on en d\'eduit que pour $x$ g\'en\'erique dans
$D$, on a
$T_x(x+K)\subset T_xD$. Le lemme 3.4 entra\^ine alors $D+K=D$. Lorsque $X$ est
simple (et en
particulier si $\dim (X)=1$), on en d\'eduit une contradiction. Sinon, l'image
de $D$ dans $X/K$ est
encore un diviseur et l'hypoth\`ese de r\'ecurrence entra\^ine encore une
contradiction. Ceci termine
la d\'emonstration du lemme.\cqfd

\ind Revenons \`a la d\'emonstration de la proposition. Notons
$Y\buildrel{h}\over{\ra}
X'\buildrel{q}\over{\ra} X$ la factorisation de Stein de $M$, o\`u $X'$ est
normal, $q$ fini
surjectif, et o\`u les fibres de $h$ sont connexes, de dimension $>0$. Si $q$
est ramifi\'e, on
choisit une composante $D'$ de son lieu de ramification; par le \theo\ de
puret\'e  ([Z]), $D=q(D')$
est un diviseur de $X$. En un point g\'en\'erique $x'$ de $D'$, le morphisme
$D'\ra D$ induit par $q$
est lisse, de sorte que l'image de la diff\'erentielle $T_{x'}q:T_{x'}X'\ra
T_{q(x')}X$ contient
$T_{q(x')}D$. Comme $T_{x'}q$ n'est pas surjective, son image est donc
$T_{q(x')}D$. Prenant pour
$Z$ une composante irr\'eductible de $h^{-1}(D')$, on obtient une contradiction
avec le lemme 4.2.
Il en r\'esulte que $q$ est non ramifi\'e, donc
\'etale ([G1], Exp. I, th. 9.5), et que $X'$ est une \va\ ([Mu1], p. 167).\cqfd

\ind Le \theo\ suivant g\'en\'eralise un r\'esultat
de Barth ([B2], Satz 2).

{\bf Th\'eor\`eme 4.3.} -- {\it Soient $X$ une \va\ simple, $V$ une vari\'et\'e
irr\'eductible,
$f:V\ra X$ un morphisme et $W$ une sous-vari\'et\'e irr\'eductible de $X$. On
suppose que $\dim \bigl(
f(V)\bigr) +2\dim (W)\ge 2\dim (X)$. Alors $f^{-1}(W)$ est connexe et, si $V$
est de plus unibranche,
le morphisme induit $\pi_1^{\rm
alg}\bigl( f^{-1}(W)\bigr)\ra\pi_1^{\rm alg}(V)$ est surjectif.}

{\bf Remarque 4.4.} La premi\`ere conclusion du \theo\ ne
subsiste pas lorsqu'on suppose seulement $\dim \bigl( f(V)\bigr) +\dim (W)>\dim
(X)$, comme le
montre l'exemple suivant d\^u \`a Barth ([B2]): soit $X'$ une \va\ simple de
dimension $7$ et $W'$ une
sous-vari\'et\'e de $X'$ lisse connexe de dimension $3$. Il existe un
sous-groupe fini $G$ non trivial
de $X'$ tel que $W'\cap W'_{\gamma}$ soit vide pour tout \'el\'ement non nul
${\gamma}$ de $G$. On pose
$X=X'/G$ et on note $W$ l'image de $W'$ dans $X$ par la projection canonique.
Alors $W$ est lisse et
il existe une sous-vari\'et\'e $V$ irr\'eductible de $X$ de dimension $5$
(intersection compl\`ete de
deux hypersurfaces), telle que $V\cap W$ soit r\'eunion de $\Card (G)$ courbes
disjointes.

{\bf D\'emonstration du \theo .} On note $\eta_V:{\tilde V}\ra V$ et
$\eta_W:{\tilde W}\ra W$ les
normalisations. Soit ${\tilde V}\times
{\tilde W}\buildrel{h}\over{\ra} X'\buildrel{q}\over{\ra} X$ la factorisation
du morphisme $(\tilde
v,\tilde w)\mapsto \bigl( f\eta_V(\tilde v)-g\eta_W(\tilde w)\bigr)$ fournie
par la proposition 4.1.
On note $G$ le noyau de $q$, on choisit un point ${\tilde v}_0$ de $\tilde V$,
un point ${\tilde
w}_0$ de $\tilde W$, et on d\'efinit des morphismes $f':{\tilde V}\ra X'$ et
$g':{\tilde W}\ra X'$ par
$f'(\tilde v)=h({\tilde v},{\tilde w}_0)$ et $g'(\tilde w)=-h({\tilde
v}_0,{\tilde w})$. L'image du
morphisme $(\tilde v,\tilde w)\mapsto h(\tilde v,\tilde w)-f'(\tilde
v)+g'(\tilde w)$ est alors
contenue dans l'ensemble fini $q^{-1}\bigl( qh({\tilde v}_0,{\tilde
w}_0)\bigr)$, donc est r\'eduite
au seul point $\alpha=-h({\tilde v}_0,{\tilde w}_0)$. On a donc: $$h(\tilde
v,\tilde w)=f'(\tilde
v)-g'(\tilde w)+\alpha\ .$$ \ind Remarquons que $f^{-1}(W)$ est l'image par
$\eta_Vpr_1$ de
$(qh)^{-1}(0)=h^{-1}(G)$. Il suffit donc de montrer que $pr_1h^{-1}(G)$ est
connexe. Pour
tout ${\gamma}\in G$, l'ensemble\break $pr_1h^{-1}({\gamma})=f'^{-1}(g'(\tilde
W)+{\gamma}-\alpha)$
est connexe. Par la remarque 3.2.2), l'intersection de $f'(\tilde V)$,
$g'(\tilde W)+{\gamma}-\alpha$
et $g'(\tilde W)-\alpha$ est non vide, de sorte que $pr_1h^{-1}({\gamma})$
rencontre $pr_1h^{-1}(0)$.
Cela entra\^ine que  $pr_1h^{-1}(G)$, donc aussi $f^{-1}(W)$, est connexe, ce
qui termine la
d\'emonstration du premier point du \theo .

\ind Supposons maintenant $V$ unibranche. Soit $\pi:V'\ra V$ un rev\^etement
\'etale connexe;
$V'$ est alors irr\'eductible (1.6). Par (1.2), il s'agit de montrer que
$f^{-1}(W)\times _VV'$ est
connexe. Or cela r\'esulte du premier point, puisque cet ensemble est
$(f\pi)^{-1}(W)$.\cqfd

\ind Le \theo\ suivant (analogue de [FL], th. 4.1) montre que, moyennant des
hypoth\`eses
suppl\'emen\-taires sur les singularit\'es des vari\'et\'es, on peut obtenir
des r\'esultats plus
forts.

{\bf Th\'eor\`eme 4.5.} -- {\it Soient $X$ une \va , $V$ et $W$ deux
vari\'et\'es irr\'eductibles et
$f:V\ra X$ et $g:W\ra X$ deux morphismes. On suppose que $V$ est unibranche,
que $f$ est minimal au
sens de (1.8), et que   $\bigl( f(V),g(W)\bigr)$ remplit strictement $X$ (\cf\
(1.9)). Alors $V\times
_XW$ est connexe. De plus, si on note $\iota$ le plongement $V\times
_XW\lhook\joinrel\mathrel{\ra}V\times W$ et $M:V\times W\ra X$ le morphisme
$(v,w)\mapsto
\bigl( f(v)-g(w)\bigr) $, alors: $$\pi_1^{\rm alg}( V\times
_XW)\buildrel{\iota_*}\over{\lra}\pi_1^{\rm alg}( V\times
W)\buildrel{M_*}\over{\lllra}
\pi_1^{\rm alg}(X)\ra 0$$ est un complexe exact en $\pi_1^{\rm alg}(X)$. Si on
suppose de plus $W$
unibranche, il est exact.}

{\bf Remarque 4.6.} La premi\`ere conclusion ne subsiste pas lorsqu'on ne fait
aucune hypoth\`ese
sur les singularit\'es de $V$, comme le montre l'exemple de Barth (remarque
4.4).

{\bf D\'emonstration du \theo .} Soient $\eta_V:{\tilde V}\ra V$ et
$\eta_W:{\tilde W}\ra W$ les
normalisations et ${\tilde V}\times {\tilde W}\buildrel{h}\over{\ra}
X'\buildrel{q}\over{\ra} X$ la
factorisation de $M\circ (\eta_V,\eta_W)$ fournie par le lemme 4.3. Par (1.8),
le morphisme
$\tilde V\ra V\ra X$ est encore minimal, de sorte que $q$ est un isomorphisme;
$V\times
_XW$ est alors l'image de $h^{-1}(0)$ par un morphisme continu, donc est
connexe. Cela montre le
premier point.

\ind Pour le second, on remarque que $V\times _XW$ est la fibre en
$0$ de $M$, de sorte qu'on a bien un complexe. Ensuite, si $X'\ra X$ est une
isog\'enie, ce qui
pr\'ec\`ede montre que $V\times _XX'$ est connexe; le morphisme $\pi_1^{\rm
alg}( V)\ra\pi_1^{\rm
alg}( X)$ est donc surjectif par (1.2). Il en est alors de m\^eme de
$\pi_1^{\rm alg}( V\times
W)\ra\pi_1^{\rm alg}( X)$.

\ind Pour montrer l'exactitude en $\pi_1^{\rm alg}( V\times W)$, on
peut supposer $V$ et $W$ normales par (1.6) et (1.8). On se donne un
rev\^etement \'etale
connexe $\pi:Y\ra V\times W$, on pose $Z=V\times _XW$ et on suppose que le
rev\^etement \'etale induit $\pi^{-1}(Z)\ra Z$ a une section.
Par [G1], Exp. V, prop. 6.11, il s'agit de montrer qu'il existe un rev\^etement
\'etale $X'\ra X$ et
un morphisme $(V\times W)\times _XX'\ra Y$ au-dessus de $V\times W$. Comme $Y$
est irr\'eductible
normale (1.6), la proposition 4.1 fournit un diagramme commutatif:
$$\diagram{Y&\hfl{h}{}&X'\cr
\vfl{\pi}{}&&\vfl{}{q}\cr
V\times W&\hfl{M}{}&X\cr}$$
o\`u les fibres de $h$ sont connexes et o\`u $q$ est une isog\'enie, dont on
note $K$
le noyau et $k$ le degr\'e. Les morphismes $\pi$ et $h$ se factorisent \`a
travers un morphisme
$\rho :Y\ra (V\times W)\times _XX'$, qui, puisque $(V\times W)\times _XX'$ est
connexe (par le premier point, puisque $V\times W$ est unibranche), est \'etale
surjectif ([G1], Exp.
V, prop. 3.5); on peut donc \'ecrire $h=h'\rho $ et $\pi=\pi'\rho $. D'autre
part,
$\pi^{-1}(Z)=h^{-1}(K)$ a $k$ composantes connexes $Y_1,\ldots ,Y_k$ et les
composantes connexes de
$\pi'^{-1}(Z)=h'^{-1}(K)$ sont $\rho (Y_1),\ldots ,\rho (Y_k)$, avec $\rho
^{-1}\bigl( \rho
(Y_i)\bigr) =Y_i$ pour tout $i$. Or par hypoth\`ese, un des rev\^etements
\'etales $\rho ^{-1}\bigl(
\rho (Y_i)\bigr)\ra \rho (Y_i)$ a une section; c'est donc un isomorphisme par
[G1], Exp. I, cor. 5.3.
Comme $V\times W$ et $\rho ^{-1}\bigl( \rho (Y_i)\bigr)$ sont connexes, cela
prouve que $\rho $ est
un isomorphisme ([G1], Exp. I, cor. 10.10) et son inverse est le morphisme
cherch\'e.\cqfd

{\bf Corollaire 4.7.} -- {\it Soient $X$ une \va\ simple, $V$ une vari\'et\'e
irr\'eductible unibranche et $f:V\ra X$ un morphisme non constant, minimal au
sens de (1.8). Alors,
le morphisme induit $\pi_1^{\rm alg}(V)\ra\pi_1^{\rm alg}(X)$ est surjectif.
Lorsque $k=\C$, il en
est de m\^eme du morphisme $\pi_1(V)\ra\pi_1(X)$.}

{\bf D\'emonstration.} La premi\`ere assertion d\'ecoule du \theo , puisque
$(V,X)$ remplit
strictement $X$ (\cf\ (1.9)). Lorsque $k=\C$, on raisonne comme dans la
d\'emonstration du cor.
3.5.\cqfd

{\bf 5. Groupe fondamental des sous-vari\'et\'es}

\ind Comme dans [FL] \pa 5, nous appliquons le \theo\ de connexit\'e \`a
l'\'etude du groupe fondamental des sous-vari\'et\'es d'une \va . Le \theo\
suivant est l'analogue
de [FL], th. 5.1.

{\bf Th\'eor\`eme 5.1.} -- {\it Soient $X$ une \va\ simple, $V$ une vari\'et\'e
irr\'eductible
unibranche et $f:V\ra X$ un morphisme non ramifi\'e. On suppose que $\dim
\bigl( f(V)\bigr) >{1\over
2}\dim (X)$. Alors il existe une \va\ $X'$, une isog\'enie $q:X'\ra X$ et un
plongement $f':V\ra X'$
tels que $f=qf'$.}

{\bf D\'emonstration.} Soit $f:V\buildrel{f'}\over{\ra}
X'\buildrel{q}\over{\ra} X$ une factorisation
de $f$ avec $f'$ minimal (1.8). On suit la d\'emonstration de [FL], p. 46: le
\theo\
4.5 implique que $V\times _{X'}V$ est connexe. Comme $f'$ est non ramifi\'e, la
diagonale $\Delta_V$
est ouverte dans $V\times _{X'}V$. Etant aussi ferm\'ee, on a
$\Delta_V=V\times _{X'}V$, de sorte
que $f'$ est bijective. On conclut en remarquant qu'un morphisme non ramifi\'e
bijectif est un
plongement ([GD2], 8.11.5 et [GD1], 17.2.6).\cqfd

\ind On en d\'eduit la g\'en\'eralisation suivante de la deuxi\`eme partie du
corollaire 2.2:

{\bf Corollaire 5.2.} -- {\it Soient $X$ une \va\ simple et $V$ une
sous-vari\'et\'e irr\'eductible
unibranche de $X$. On suppose que $\dim
(V)>{1\over 2}\dim (X)$. Alors, le morphisme induit\break $\pi_1^{\rm
alg}(V)\ra\pi_1^{\rm
alg}(X)$ est un isomorphisme.}

{\bf D\'emonstration.} Par la remarque 3.6.2), l'inclusion
$V\lhook\joinrel\mathrel{\ra}X$
est minimale. Soit $\pi:V'\ra V$ un rev\^etement \'etale connexe; $V$ \'etant
irr\'eductible et
unibranche, il en est de m\^eme de $V'$ (1.6). Le \theo\ pr\'ec\'edent
s'applique au morphisme non
ramifi\'e
 $V'\buildrel{\pi}\over{\ra}
V\buildrel{\iota}\over{\lhook\joinrel\mathrel{\ra}} X$, qui se factorise
donc en:
$$\diagram{V'&\buildrel{\iota'}\over{\lhook\joinrel\mathrel{\hfl{}{}}}&X'\cr
\vfl{\pi}{}&&\vfl{}{q}\cr
V&\buildrel{\iota}\over{\lhook\joinrel\mathrel{\hfl{}{}}}&X\cr}$$
o\`u $\iota'$ est un plongement et $q$ une isog\'enie. Le morphisme induit
$V'\ra V\times
_XX'$ est alors un plongement \'etale. Comme $V\times _XX'$ est connexe (\theo\
4.5), c'est un
isomorphisme ([G1], Exp. I, cor. 5.2). Par (1.3), ceci montre l'injectivit\'e
du morphisme
$\pi_1^{\rm alg}(V)\ra\pi_1^{\rm alg}(X)$; celui-ci est d'autre part surjectif
par le corollaire
3.5.\cqfd

{\bf Remarques 5.3.} 1) J'ignore si la restriction sur les singularit\'es de
$V$ est n\'ecessaire.

\ind 2) Lorsque $V$ est de plus normal, il ressort du corollaire et de (1.5)
que le
morphisme de restriction $\Pic ^0(X)\ra\Pic ^0(V)$ est un isomorphisme. Si $V$
est lisse, le
corollaire 2.2 montre que la restriction $\Pic (X)\ra\Pic (V)$ est {\it
injective}. Il serait
int\'eressant de savoir si cette conclusion subsiste si l'on permet \`a $V$
d'\^etre singulier.

{\bf 6. Le probl\`eme de Zariski pour les \vas }

\ind Soient $X$ une \va\ et $D$ un diviseur de $X$. Le probl\`eme de Zariski
est l'\'etude du
groupe fondamental alg\'ebrique $X- D$. En suivant la d\'emarche de [F2], on va
voir qu'il est
possible de retrouver certains cas particuliers des r\'esultats tr\`es
g\'en\'eraux de [N].

\ind On dira que le
diviseur $D$ de $X$ est {\it \`a croisements normaux en codimension $1$} s'il
existe un sous-sch\'ema
ferm\'e $Z$ de $D$ de codimension $2$ dans $D$ tel que $D- Z$ soit un diviseur
\`a
croisements normaux dans $X- Z$.

{\bf Th\'eor\`eme 6.1.} -- {\it Soient $X$ une \va\ de dimension $>1$ et $D$ un
diviseur de $X$ \`a
croisements normaux en codimension $1$. On suppose qu'aucune composante de $D$
n'est une sous-\va\
de $X$. Alors le noyau du morphisme canonique  $\pp (X- D)\ra\pi_1^{\rm
alg}(X)$ est
ab\'elien.}

{\bf Remarque 6.2.} Il est clair que le r\'esultat ne subsiste pas en
g\'en\'eral si
$D$ contient une sous-\va\ $K$ de $X$. Notons  $E$ la courbe elliptique $X/K$;
on a alors des
surjections: $$\Ker \bigl(\pi_1^{\rm alg} (X- D)\!\ra\!\pi_1^{\rm
alg}(X)\bigr)\dra\Ker\bigl(\pi_1^{\rm alg} (X- K)\!\ra\!\pi_1^{\rm
alg}(X)\bigr)
\dra\Ker\bigl(\pi_1^{\rm alg} (E-
\{ 0\} )\!\ra\!\pi_1^{\rm alg}(E)\bigr)$$\ind  Le groupe fondamental
alg\'ebrique de $E-
\{ 0\}$ est isomorphe au compl\'et\'e profini de $\Z \mathrel*\Z$, dont le
noyau du morphisme canonique vers $\pi_1^{\rm alg}(E)\isom \hat{\Z}$ n'est pas
ab\'elien.

{\bf D\'emonstration du \theo .} Soient $V$ une vari\'et\'e irr\'eductible
normale et $f:V\ra X$ un
rev\^etement Galoisien de groupe $G$, dont le discriminant
est contenu dans $D$. Il
existe une factorisation  $f:V\buildrel{f'}\over{\ra} X'\buildrel{q}\over{\ra}
X$, o\`u $q$ est une
isog\'enie et o\`u $f'$ est un rev\^etement Galoisien de groupe $G'\subset G$,
minimal au
sens de (1.8). Le diviseur $D'=q^{-1}(D)$ de $X'$ est  encore \`a croisements
normaux
en codimension $1$. Soit $B$ une composante irr\'eductible de $D'$. Nous allons
montrer que
$f'^{-1}(B)$ est irr\'eductible.

\ind Soit $\tilde B\ra B$ la normalisation. On v\'erifie \`a l'aide du
lemme d'Abhyankar ([G1], 3.6, Exp. X) que $(\tilde B\times _{X'}V)_{\rm red}$
est non singulier en
codimension $1$. Comme $V$ est normale, que $f'$ est minimale et que $B$
engendre $X'$, il r\'esulte
du \theo\ 4.5 que $\tilde B\times _{X'}V$ est connexe. Comme il est de plus de
Cohen-Macaulay,
il s'ensuit que ce sch\'ema est irr\'eductible ([H1]). Il en de m\^eme pour sa
projection $f'^{-1}(B)$
sur $V$.

\ind Les composantes de $D'$ se rencontrant deux \`a deux, il en est de m\^eme
pour leurs images inverses dans $V$. Le lemme d'Abhyankar entra\^ine que les
inerties
correspondantes commutent deux \`a deux ([S]). Elles engendrent donc un
sous-groupe ab\'elien $I$
de $G'$. Posons $X''=V/I$. Le morphisme $f'$ se factorise en $V\ra
X''\buildrel{f''}\over{\ra} X'$ et
$f''$ est alors non ramifi\'e en codimension $1$, donc non ramifi\'e d'apr\`es
le \theo\ de puret\'e.

\ind Posons $D''=f''^{-1}(D')$. On a un diagramme commutatif:
$$\diagram{\pp (X- D)&\ldra&G\cr
\cup&&\cup\cr
\pp (X''- D'')&\ldra&I\cr}$$
\ind D'autre part, $\pp (X''- D'')$ contient le noyau $N$ du morphisme
canonique\break
$\pp (X- D)\ra\pi_1^{\rm alg}(X)$ ([G1], cor. 6.7, Exp. V). Il s'ensuit que
l'image de $N$ dans $G$
est contenue dans $I$, donc est ab\'elienne.

\ind Nous avons donc montr\'e que pour tout sous-groupe ouvert distingu\'e
$\Gamma$ de\break $\pp
(X- D)$, le groupe $N/(N\cap\Gamma)$ est ab\'elien. Comme $N$ est ferm\'e, il
s'identifie \`a
$\limproj_ {\,\Gamma}N/(N\cap\Gamma)$  ([Bo], III, \pa 7, n\up{o}2 , prop. 3).
C'est donc un
groupe ab\'elien.\cqfd

\ind Lorsque $k=\C $, Nori montre dans [N] que le noyau du
morphisme $\pi_1(X- D)\ra \pi_1 (X)$ entre groupes fondamentaux {\it
topologiques} est ab\'elien de
type fini, et que son centralisateur est d'indice fini. On peut pr\'eciser son
r\'esultat de la fa\c
con suivante.

{\bf Proposition 6.3.} -- {\it Sous les hypoth\`eses du \theo , et lorsque
$k=\C$, le noyau $N$ du
morphisme canonique  $\pi_1(X- D)\ra\pi_1(X)$ est ab\'elien libre de type
fini.}

{\bf Remarques 6.4.} 1) Soient $\tilde D_1,\ldots ,\tilde D_r$ les composantes
irr\'eductibles de la
normalisation de $D$, et $K_i$ le noyau du morphisme $\Pic ^0(X)\ra\Pic
^0(\tilde D_i)$. Alors le rang
de $N$ est $\sum_{i=1}^r\Card (K_i)$.

\ind 2) Lorsque les composantes du diviseur $D$ sont de plus normales, le rang
de $N$ est le nombre
$r$ de composantes irr\'eductibles de $D$, et Nori montre que l'extension:
$$0\ra \Z ^r \ra \pi_1(X- D)\ra\pi_1(X)\ra 0\ ,$$est {\it centrale} et que sa
classe correspond au
$r$--uplet des classes fondamentales des composantes de $D$ via l'isomorphisme
$H^2\bigl( \pi_1(X),\Z
^r\bigr)\isom H^2(X,\Z )^r$.

{\bf D\'emonstration de la proposition.} Soit $k$ le nombre maximum de
composantes
irr\'eductibles de l'image inverse de $D$ par une isog\'enie $X'\ra X$. On
v\'erifie que ce nombre
est fini et qu'il se calcule comme dans la remarque ci-dessus. Les inerties
\'etant cycliques, la
d\'emonstration du \theo\ montre que tout quotient fini de $N$ (donc aussi $N$)
peut \^etre
engendr\'e par $k$ \'el\'ements. Soit $q:X'\ra X$ une isog\'enie telle que
$D'=q^{-1}(D)$ ait $k$
composantes. On v\'erifie que $N$ est isomorphe au noyau $N'$ du morphisme
$\pi_1(X'-D')\ra\pi_1(X')$. D'autre part, il existe une
surjection de $N'$ sur le noyau de $H_1(X'-D',\Z )\ra H_1(X',\Z )$. Un argument
classique de dualit\'e
montre que ce dernier est isomorphe au conoyau du morphisme $H_2(X',\Z )\ra
H^{2\dim (D')}(D',\Z
)\isom\Z ^k$, qui envoie $\gamma$ sur $\bigl( d_1(\gamma),\ldots
,d_k(\gamma)\bigr)$, o\`u
$d_1,\ldots ,d_k$ sont les classes fondamentales des composantes de $D'$. Soit
$M:X'\ra X'$ la
multiplication par un entier $m$. Posons $D''=M^{-1}(D')$. De nouveau, $N'$ est
isomorphe au noyau
$N''$ du morphisme $\pi_1(X''-D'')\ra\pi_1(X'')$ et les composantes de $D''$
sont les images inverses
de celles de $D'$, donc leurs classes sont divisibles par $m^2$. Il ressort de
ce qui
pr\'ec\`ede qu'il existe une surjection de $N''$ sur $\Z ^k/(m^2,\ldots ,m^2)$.
Ceci \'etant vrai
pour tout $m$, il existe une surjection de  $N$ sur $\Z ^k$. Comme $N$ est
engendr\'e par $k$
\'el\'ements, il est isomorphe \`a $\Z ^k$.\cqfd

{\bf 7. Rev\^etements des sous-vari\'et\'es des \vas }

\ind Suivant toujours [FL], on s'int\'eresse maintenant aux vari\'et\'es qui
sont
rev\^etements finis ramifi\'es de petit degr\'e d'une sous-vari\'et\'e d'une
\va\ simple. Pour
\'enoncer notre r\'esultat, dont la d\'emonstration est calqu\'ee sur celle du
\theo\ principal de
[GL], nous utiliserons la notion de {\it degr\'e local} d'un morphisme fini
$f:V\ra W$ en un point
$v\in V$, telle qu'elle est d\'efinie dans [GL] et [Mu2]. Intuitivement, c'est
le nombre de feuillets
de $f$ qui se rejoignent en $v$. On le note $e_f(v)$.

{\bf Th\'eor\`eme 7.1.} -- {\it Soient $X$ une \va\ simple, $W$ une
sous-vari\'et\'e lisse de $X$,
$V$ une vari\'et\'e irr\'eductible normale et $f:V\ra W$ un rev\^etement fini
de degr\'e $d$. On
suppose le morphisme compos\'e $V\ra W\lhook\joinrel\mathrel{\ra}X$ est minimal
au sens de
(1.8). Alors il existe un point $v\in V$ tel que: $$e_f(v)\ge \min \bigl(
d,2\dim
(W)-\dim (X)+1\bigr)\ .$$}
{\bf D\'emonstration.} Pour tout entier $l$,
l'ensemble: $$R_l=\{ v\in V \bigm| e_f(v)>l\}\ ,$$ est ferm\'e dans $V$ ([GL],
lemma 1). De plus,
comme $W$ est lisse, le \theo\ 2.2 de [L2] (\cf\ aussi [GL], p. 58), entra\^ine
que $R_l$ est soit
vide, soit partout de codimension $\le l$ dans $V$. Il s'agit de montrer que
$R_l$ est non vide pour
$l\le \min \bigl( d-1,2\dim (W)-\dim (X)\bigr)$. Nous proc\'edons par
r\'ecurrence sur $l$. On a
$R_0=V$; on suppose que $R_{l-1}$ est non vide, et on en choisit une composante
irr\'eductible $R$.
La projection $V\times _X R\ra R$ est finie, de sorte que la diagonale
$\Delta_R$ est une composante
irr\'eductible de $V\times _X R$. Comme:
$$\dim (R)\ge \dim (V)-(l-1)\ge \dim (V)-2\dim (W)+\dim(X)+1>\dim(X)-\dim (V)\
,$$
le \theo\ 4.5 entra\^ine que $V\times _X R$ est connexe. Si $V\times _X
R=\Delta_R$, alors $R\subset
R_{d-1}\subset R_l$, ce qui permet de conclure. Sinon, il existe une composante
irr\'eductible $T$ de
$V\times _X R$ distincte de $\Delta_R$, qui rencontre la diagonale. Comme dans
[GL], p. 57, on
conclut que pour tout point $(v,v)$ de $T\cap \Delta_R$, on a $e_f(v)> l$,
\cad\ $v\in R_l$.\cqfd

\ind Comme dans [GL], on en d\'eduit:

{\bf Corollaire 7.2.} -- {\it Soient $X$ une \va\ simple, $W$ une
sous-vari\'et\'e lisse de $X$,
$V$ une vari\'et\'e irr\'eductible unibranche et $f:V\ra W$ un rev\^etement
fini de degr\'e\break
$d\le 2\dim (W)-\dim (X)$. Alors, le
morphisme induit $\pi_1^{\rm alg}(V)\ra\pi_1^{\rm alg}(X)$ est injectif et son
image est isomorphe
\`a $\pi_1^{\rm alg}(X)$.}

{\bf D\'emonstration.} Par (1.6), on peut supposer $V$ normale. Toute
isog\'enie $X'\ra X$
induit une injection $\pi_1^{\rm alg}(X')\ra\pi_1^{\rm alg}(X)$ (1.4) dont
l'image est isomorphe \`a
$\pi_1^{\rm alg}(X)$ (1.7). On peut donc supposer que le morphisme compos\'e
$g:V\ra
W\lhook\joinrel\mathrel{\ra}X$ est minimal. Soit $\pi:V'\ra V$ un rev\^etement
\'etale connexe; $V$
\'etant irr\'eductible et normale, il en est de m\^eme de $V'$ (1.6).
Soit $V'\buildrel{g'}\over{\ra} X'\buildrel{p}\over{\ra} X$ une factorisation
de $g\pi$ avec $g'$
minimal (1.8). Posons $W'=g'(V')$, et notons $f':V'\ra W'$ le rev\^etement
induit et $d'$ son
degr\'e. Le \theo\ s'applique \`a $f'$ et il existe $v'\in V'$ tel que:
$$e_{f'}(v')\ge \min \bigl(
d',2\dim (W')-\dim (X')+1\bigr)\ .$$Comme $p$ et $\pi$ sont \'etales et que $d$
est le degr\'e de
$f$, on a: $$e_{f'}(v')=e_{pf'}(v')=e_{f\pi}(v')=e_f\bigl( \pi(v')\bigr)\le d\
.$$ Comme $d\le2\dim
(W)-\dim (X)$, il s'ensuit que $d\ge d'$, donc que:
$$\deg (\pi)\le \deg (p_{|W'})\le \deg (p)\ .$$ Or $\pi$ se
factorise \`a travers le morphisme \'etale $V\times _XX'\ra V$, qui est de
m\^eme degr\'e que $p$.
Comme $V\times _XX'$ est connexe (\theo\ 4.5) donc irr\'eductible (1.6), le
morphisme induit $V'\ra
V\times _XX'$ est un isomorphisme. Par (1.3), cela d\'emontre le
corollaire.\cqfd

{\bf Remarques 7.3.} 1) La conclusion du corollaire ne subsiste en g\'en\'eral
pas sans hypoth\`ese
sur les singularit\'es de $V$, comme le montre la construction de [FL], note
(2), p. 56.

\ind 2) On a besoin de la lissit\'e de $W$  uniquement pour
pouvoir appliquer le r\'esultat de Lazarsfeld. Plus pr\'ecis\'ement, il faut
savoir que les lieux
$R_l$ sont de codimension $\le l$ pour $l<2\dim (W)-\dim (X)<d$. Un \theo\ de
puret\'e de
Grothendieck ([G3], Exp. X, th. 3.4) dit que c'est le cas pour $V$ normal, $W$
localement intersection
compl\`ete, $l=1$ et $\codim\bigl(\Sing (W)\bigr)\ge 3$. Cela sugg\`ere que le
corollaire devrait
s'\'etendre au cas  $W$ localement intersection compl\`ete et $\dim\bigl(\Sing
(W)\bigr)< \codim (W)$.

{\bf Corollaire 7.4.} -- {\it Soient $X$ une \va\ simple, $W$ une
sous-vari\'et\'e lisse de $X$,
$V$ une vari\'et\'e irr\'eductible et $f:V\ra W$ un rev\^etement fini de
degr\'e $d$. Si $V_1$ et
$V_2$ sont deux sous-vari\'et\'es de $V$ telles que:$$\dim (V_1)+\dim (V_2)\ge
2\dim (X)-\dim
(V)+d-1\ ,$$alors $V_1\cap V_2\ne\emptyset$. En particulier, tout morphisme de
$V$ dans une
vari\'et\'e de dimension $\le 3\dim (V)-2\dim (X)-d$ est constant.}

{\bf D\'emonstration du corollaire.}  La d\'emonstration suit [L1], remarque
2.3. On peut
supposer $V$ normale et $d\le
2\dim (W)-\dim (X)$. Par le \theo\ 7.1 et le r\'esultat de [L2] utilis\'e plus
haut, il existe une
sous-vari\'et\'e $R$ de $V$ de codimension $\le d-1$ telle que $f$ soit
bijective au-dessus de
$f(R)$. La remarque 3.2.2) entra\^ine alors que $f(V_1)$, $f(V_2)$ et $f(R)$ se
rencontrent, donc
aussi $V_1$, $V_2$ et $R$. La derni\`ere assertion en r\'esulte
imm\'ediatement: si $V\ra T$ est
surjectif et si $\dim (T)>0$, il suffit de prendre pour $V_1$ l'image inverse
d'un diviseur et pour
$V_2$ l'image inverse d'un point hors de ce diviseur.\cqfd

\ind La d\'emonstration du \theo\ suivant suit celle d'un r\'esultat analogue
de
Lazarsfeld ([F3], th., p. 151).

{\bf Th\'eor\`eme 7.5.} -- {\it Soient $X$ une \va\ simple, $V$ une vari\'et\'e
irr\'eductible
unibranche et $f:V\ra X$ un morphisme. On suppose qu'il existe une
sous-vari\'et\'e $Z$ de $f(V)$
v\'erifiant $\dim (Z)>\codim \bigl( f(V)\bigr) $ telle que le morphisme induit
$f^{-1}(Z)\ra Z$ soit
une bijection ensembliste. Alors, le
morphisme induit $\pi_1^{\rm alg}(V)\ra\pi_1^{\rm alg}(X)$ est injectif et son
image est isomorphe
\`a $\pi_1^{\rm alg}(X)$.}

{\bf Remarque 7.6.} Le \theo\ s'applique en particulier aux
rev\^etements cycliques des sous-vari\'et\'es de $X$ de dimension $\ge1+{1\over
2}\dim (X)$.

{\bf D\'emonstration.} On peut comme d'habitude supposer $f$ minimal, et $Z$
int\`egre. Soit $\tilde
Z$ la normalisation de $Z$; on pose $Z'=\tilde Z\times _XV$. Nos hypoth\`eses
entra\^inent que le
morphisme induit $Z'_{\rm red}\ra\tilde Z$ est fini et birationnel; c'est donc
un isomorphisme par le
\theo\ principal de Zariski ([GD2], 8.12.10.1). Le morphisme  $\pi_1^{\rm
alg}(Z'_{\rm
red})\ra\pi_1^{\rm alg}(Z')$ \'etant surjectif (1.2), l'intersection de l'image
du morphisme:
$$\pi_1^{\rm alg}(Z')\ra\pi_1^{\rm alg}(\tilde Z\times V)\isom\pi_1^{\rm
alg}(\tilde Z)\times
\pi_1^{\rm alg}(V)$$avec $\{ 1\}\times\pi_1^{\rm alg}(V)$ est r\'eduite \`a $\{
(1,1)\}$. Comme
cette image est le noyau du morphisme canonique $\pi_1^{\rm alg}(\tilde Z\times
V)\ra\pi_1^{\rm alg}(X)$ (\theo\ 4.5), il s'ensuit que le morphisme $\pi_1^{\rm
alg}(V)\ra\pi_1^{\rm
alg}(X)$ est injectif. Il est d'autre part surjectif par le corollaire 4.7, ce
qui conclut la
d\'emonstration.\cqfd

\ind On termine ce chapitre avec la d\'emonstration d'un r\'esultat analogue
\`a un \theo\ de Ein
sur les rev\^etements des espaces projectifs. Elle utilise des r\'esultats
puissants de
 Kawamata, Koll\'ar et Mori, mais je pense qu'il devrait exister une
d\'emonstration
commune (plus simple) \`a base de connexit\'e.

{\bf Th\'eor\`eme 7.7.} -- {\it Soient $X$ est une \va\ simple, $V$ une
vari\'et\'e lisse  et $f:V\ra X$
un rev\^etement fini ramifi\'e. Alors la ramification de $f$, \cad\ le
faisceau canonique de $V$, est ample.}

{\bf Remarque 7.8.} Le \theo\ de connexit\'e 4.5 entra\^ine que la ramification
de $f$ rencontre
toute courbe de $V$.

{\bf D\'emonstration.} Par [K1], theorem 13, il existe une sous-vari\'et\'e
ab\'elienne $Y$ de $X$,
un morphisme fini $W\ra X/Y$ avec $\kappa (V)=\dim (W)$,
 et des rev\^etements \'etales
${\tilde Y}\ra Y$ et ${\tilde Y}\times W\ra V$. Comme $X$ est simple, on a soit
$Y=X$, auquel
cas $f$ est \'etale, ce qui contredit l'hypoth\`ese; soit $Y=0$, et $\kappa
(V)=\dim (V)$, de
sorte que $V$ est de type g\'en\'eral. Le
 \theo\ d\'ecoule alors du lemme suivant.\cqfd

{\bf Lemme 7.9.} -- {\it Soient $V$ une vari\'et\'e projective lisse de type
g\'en\'eral. Alors
soit $K_V$ est ample, soit $V$ contient une courbe rationnelle.}

{\bf D\'emonstration.} On suit [Ko]. Par [M], theorem 1.4, soit $V$ contient
une courbe
rationnelle, soit $K_V$ est nef. Dans le second cas, il d\'ecoule de [CKM],
(9.3), que le syst\`eme
lin\'eaire $|mK_V|$
 est sans point base pour $m\gg 0$, donc d\'efinit un morphisme $\phi$
g\'en\'eriquement fini sur
son image. Si $\phi$ est fini, $K_V$ est ample; sinon, il contracte une courbe
 $C$, qui v\'erifie alors $C\cdot K_V=0$. Soit $H$ un diviseur tr\`es ample sur
$V$; son image
$\phi (H)$ est contenue dans une hypersurface de degr\'e $r$, de sorte qu'il
existe un diviseur
effectif $E$  tel que $rmK_V\equiv H+E$. On a donc $C\cdot E<0$. Il existe
alors $\epsilon$
rationnel positif tel que le diviseur $K_V+\epsilon E$ soit {\it log-canonique}
([KMM], lemma
0-2-15). Comme $(K_V+\epsilon E)\cdot C <0$, il n'est pas nef, et le \theo\ du
c\^one (\loc ,
theorem 4-2-1) entra\^ine qu'il existe une raie extr\^emale et un morphisme de
contraction (\loc ,
theorem 3-2-1) $g:V\ra W$, auquel on peut appliquer [K2], theorem 1: le lieu
des points o\`u $g$
n'est pas un isomorphisme est recouvert par des courbes rationnelles.\cqfd

{\bf 8. Conjecture}

\ind Une d\'emonstration de la conjecture suivante permettrait d'utiliser les
r\'esultats de
Lazarsfeld ([L1], th. 2.1, prop. 3.1), et entra\^inerait en particulier notre
\theo\ 7.1 sur les
indices de ramification.

{\bf Conjecture 8.1.} {\it Soient $X$ est une \va\ simple, $V$ une vari\'et\'e
lisse et $f:V\ra X$
un rev\^etement fini, minimal au sens de (1.8). On d\'efinit un faisceau
localement
libre $E$ sur $X$ comme le dual du noyau de la trace $\Tr _{V/X}:f_*{\cal
O}_V\lra{\cal O}_X$. Alors
$E$ est ample.}

\ind La r\'eponse est affirmative si $X$ est une courbe elliptique. En effet,
par un \theo\
de Hartshorne ([H2]), il suffit de d\'emontrer que tout faisceau inversible $L$
quotient de $E$ est de
degr\'e $>0$. Soit $L$ un tel faisceau; on a alors: $$0<h^0(X,E^{\vee}\otimes
L)=h^0(X,f_*{\cal
O}_V\otimes L)-h^0(X,L)=h^0(V,f^*L)-h^0(X,L)\ .$$Cela entra\^ine en particulier
$h^0(V,f^*L)>0$. Si
$\deg (L)\le 0$, alors $\deg
(f^*L)\le 0$ et $f^*L$ est donc trivial. On a alors $h^0(V,f^*L)=1$ et
$h^0(X,L)=0$, donc
$L\not\isom {\cal O}_X$, ce qui contredit l'injectivit\'e de $f^*$. On a donc
d\'emontr\'e que
$\deg (L)>0$.

\ind On peut voir cette conjecture comme une forme forte d'un \theo\ de
connexit\'e: supposons $E$
ample; si $W$ est une vari\'et\'e irr\'eductible de dimension $>0$ et si
$g:W\ra X$ est un morphisme
fini, alors $g^*E$ est ample, donc $H^0(W, g^*E^{\vee})=0$. Si on pose
$W'=V\times
_XW$, cela se traduit par $h^0(W',{\cal O}_{W'})=h^0(W,{\cal O}_W)=1$, qui est
\'evidemment plus fort
que la connexit\'e de $W'$ (assur\'ee par le \theo\ 4.5).

\saut\saut \centerline {\pc REFERENCES}\saut
\hangindent=1cm
[B1]	Barth, W., {\it Fortsetzung meromorpher Funktionen in Tori und
Komplexprojektiven R\"aumen},
Invent. Math {\bf 5} (1968), 42-62.

\hangindent=1cm
[B2]	Barth, W., {\it Verallgemeinerung des bertinischen Theorems in abelschen
Mannigfaltig\-kei\-ten}, Ann. Sc. Norm. Sup. Pisa, Serie IV, {\bf 23} (1969),
317--330.

\hangindent=1cm
[Bo]	Bourbaki,N., Topologie G\'en\'erale, Hermann, 1960.

\hangindent=1cm
[CKM] Clemens, H., Koll\'ar, J., Mori, S., {\it  Higher Dimensional Complex
Geometry\/},
Ast\'erisque 166.

\hangindent=1cm
[F1] Fulton, W., On the Topology of Algebraic Varieties, in {\it Algebraic
Geometry
Bowdoin 1985\/}, Proceedings of Symposia in Pure Mathematics {\bf 46}, Part 1,
1987.

\hangindent=1cm
[F2] Fulton, W., {\it On the Fundamental Group of the Complement of a Node
Curve\/}, Ann. Math.
 {\bf 11} (1980), 407--409.

\hangindent=1cm
[F3] Fulton, W., On Nodal Curves, in {\it Algebraic Geometry -- Open
Problems\/},
Proceedings, Ravello 1982, Springer Lecture Notes 977.

\hangindent=1cm
[FH] Fulton, W., Hansen, J. {\it A Connectedness Theorem for Projective
Varieties, with
Applications to Intersections and Singularities of Mappings\/}, Ann. of Math.
{\bf 110} (1979),
159--166.

\hangindent=1cm
[FL] Fulton, W., Lazarsfeld, R., Connectivity and its Applications in Algebraic
Geometry, in
{\it Algebraic Geometry\/}, Proceedings of the Midwest Algebraic Geometry
Conference, Chicago 1980,
Springer Lecture Notes 862.

\hangindent=1cm
[GL] Gaffney, T., Lazarsfeld, R., {\it On the Ramification of Branched
Coverings of $\P ^n$\/},
Invent. Math. {\bf 59} (1980), 53--58.

\hangindent=1cm
[G1] Grothendieck, A., {\it Rev\^etements Etales et Groupe Fondamental\/},
S.G.A. 1, Springer Lecture
Notes 224.

\hangindent=1cm
[G2] Grothendieck, A., {\it Fondements de la G\'eom\'etrie Alg\'ebrique\/},
S\'eminaire Bourbaki
1957--1962.

\hangindent=1cm
[G3] Grothendieck, A., {\it Cohomologie locale des faisceaux coh\'erents et
th\'eor\`emes de Lefschetz
locaux et globaux\/}, S.G.A. 2, North Holland, Amsterdam, 1968.

\hangindent=1cm
[GD1] Grothendieck, A., Dieudonn\'e, J., {\it El\'ements de G\'eom\'etrie
Alg\'ebrique IV, 4\/},
Publ. Math. I.H.E.S. 32, 1967.

\hangindent=1cm
[GD2] Grothendieck, A., Dieudonn\'e, J., {\it El\'ements de G\'eom\'etrie
Alg\'ebrique IV, 3\/},
Publ. Math. I.H.E.S. 28, 1966.

\hangindent=1cm
[H1]	Hartshorne, R., {\it Complete Intersections and Connectedness}, Amer. J.
Math. {\bf
84} (1962), 497--508.

\hangindent=1cm
[H2]	Hartshorne, R., {\it Ample vector bundles on curves}, Nagoya Math. J. {\bf
43}
(1971), 73--89.

\hangindent=1cm
[KMM]	Kawamata, Y., Matsuda, K., Matsuki, K., Introduction to the Minimal Model
Problem, in  {\it
Algebraic Geometry Sendai 1985}, Oda, T. editor, Adv. Stud. Pure Math. {\bf
10}, Tokyo, 1987.

\hangindent=1cm
[K1]	Kawamata, Y., {\it Characterization of Abelian Varieties}, Comp. Math.
{\bf 43}
(1981), 253--276.

\hangindent=1cm
[K2]	Kawamata, Y., {\it On the Length of an Extremal Rational Curve\/}, Inv.
Math. {\bf 105}
(1991), 609--611.

\hangindent=1cm
[Ko] Koll\'ar, J., {\it Shafarevich Maps and Automorphic Forms\/}, preprint.

\hangindent=1cm
[L1] Lazarsfeld, R., {\it A Barth-Type Theorem for Branched Coverings of
Projective Space\/}, Math.
Ann. {\bf 249} (1980), 153--162.

\hangindent=1cm
[L2] Lazarsfeld, R., Ph.D. thesis, Brown University, June 1980.

\hangindent=1cm
[M]	Mori, S., {\it Threefolds whose Canonical Bundles are not Numerically
Effective\/}, Ann.
Math. {\bf 116} (1982), 133--176.

\hangindent=1cm
[Mu1] Mumford, D., {\it Abelian Varieties\/}, Oxford University Press, 1974.

\hangindent=1cm
[Mu2] Mumford, D., {\it Algebraic Geometry I. Complex Projective Varieties\/},
Grundlehren der
mathematischen Wissenschaften 221, Springer Verlag, 1976.

\hangindent=1cm
[N]	Nori, M., {\it Zariski's conjecture and related problems\/}, Ann. Sci.
Ecole Norm. Sup. {\bf 16}
(1983), 305--344.

\hangindent=1cm
[S] Serre, J.-P., {\it Rev\^etements ramifi\'es du plan projectif (d'apr\`es S.
Abhyankar)}, S\'eminaire Bourbaki, Exp. 204, 1960.

\hangindent=1cm
[So]	Sommese, A., {\it Complex Subspaces of Homogeneous Complex Manifolds II.
Homotopy Results\/},
Nagoya Math. J. {\bf 86} (1982), 101--129.

\hangindent=1cm
[Z]	Zariski, O., {\it On the Purity of the Branch Locus of Algebraic
Functions\/},
Proc. Nat. Acad. Sci. {\bf 44} (1958), 791--796.

\bye